\documentclass[12pt, 14paper,reqno]{amsart}
\usepackage{amsmath,amsfonts,amssymb}
\usepackage{longtable}

\usepackage[mathscr]{eucal}

\usepackage{amssymb, amsmath, amsthm}
\usepackage[breaklinks]{hyperref}

\usepackage{graphicx}

\input xypic
\xyoption{all}

\makeindex
\makeglossary

\begin{document}

\baselineskip = 16pt

\newcommand \ZZ {{\mathbb Z}}
\newcommand \NN {{\mathbb N}}
\newcommand \RR {{\mathbb R}}
\newcommand \PR {{\mathbb P}}
\newcommand \AF {{\mathbb A}}
\newcommand \GG {{\mathbb G}}
\newcommand \QQ {{\mathbb Q}}
\newcommand \CC {{\mathbb C}}
\newcommand \bcA {{\mathscr A}}
\newcommand \bcC {{\mathscr C}}
\newcommand \bcD {{\mathscr D}}
\newcommand \bcE {{\mathscr E}}
\newcommand \bcF {{\mathscr F}}
\newcommand \bcG {{\mathscr G}}
\newcommand \bcH {{\mathscr H}}
\newcommand \bcM {{\mathscr M}}
\newcommand \bcI {{\mathscr I}}
\newcommand \bcJ {{\mathscr J}}
\newcommand \bcK {{\mathscr K}}
\newcommand \bcL {{\mathscr L}}
\newcommand \bcO {{\mathscr O}}
\newcommand \bcP {{\mathscr P}}
\newcommand \bcQ {{\mathscr Q}}
\newcommand \bcR {{\mathscr R}}
\newcommand \bcS {{\mathscr S}}
\newcommand \bcT {{\mathscr T}}
\newcommand \bcV {{\mathscr V}}
\newcommand \bcU {{\mathscr U}}
\newcommand \bcW {{\mathscr W}}
\newcommand \bcX {{\mathscr X}}
\newcommand \bcY {{\mathscr Y}}
\newcommand \bcZ {{\mathscr Z}}
\newcommand \goa {{\mathfrak a}}
\newcommand \gob {{\mathfrak b}}
\newcommand \goc {{\mathfrak c}}
\newcommand \gom {{\mathfrak m}}
\newcommand \gon {{\mathfrak n}}
\newcommand \gop {{\mathfrak p}}
\newcommand \goq {{\mathfrak q}}
\newcommand \goQ {{\mathfrak Q}}
\newcommand \goP {{\mathfrak P}}
\newcommand \goM {{\mathfrak M}}
\newcommand \goN {{\mathfrak N}}
\newcommand \uno {{\mathbbm 1}}
\newcommand \Le {{\mathbbm L}}
\newcommand \Spec {{\rm {Spec}}}
\newcommand \Gr {{\rm {Gr}}}
\newcommand \Pic {{\rm {Pic}}}
\newcommand \Jac {{{J}}}
\newcommand \Alb {{\rm {Alb}}}
\newcommand \alb {{\rm {alb}}}
\newcommand \Corr {{Corr}}
\newcommand \Chow {{\mathscr C}}
\newcommand \Sym {{\rm {Sym}}}
\newcommand \Prym {{\rm {Prym}}}
\newcommand \cha {{\rm {char}}}
\newcommand \eff {{\rm {eff}}}
\newcommand \tr {{\rm {tr}}}
\newcommand \Tr {{\rm {Tr}}}
\newcommand \pr {{\rm {pr}}}
\newcommand \ev {{\it {ev}}}
\newcommand \cl {{\rm {cl}}}
\newcommand \interior {{\rm {Int}}}
\newcommand \sep {{\rm {sep}}}
\newcommand \td {{\rm {tdeg}}}
\newcommand \alg {{\rm {alg}}}
\newcommand \im {{\rm im}}
\newcommand \gr {{\rm {gr}}}
\newcommand \op {{\rm op}}
\newcommand \Hom {{\rm Hom}}
\newcommand \Hilb {{\rm Hilb}}
\newcommand \Sch {{\mathscr S\! }{\it ch}}
\newcommand \cHilb {{\mathscr H\! }{\it ilb}}
\newcommand \cHom {{\mathscr H\! }{\it om}}
\newcommand \colim {{{\rm colim}\, }}
\newcommand \End {{\rm {End}}}
\newcommand \coker {{\rm {coker}}}
\newcommand \id {{\rm {id}}}
\newcommand \van {{\rm {van}}}
\newcommand \spc {{\rm {sp}}}
\newcommand \Ob {{\rm Ob}}
\newcommand \Aut {{\rm Aut}}
\newcommand \cor {{\rm {cor}}}
\newcommand \Cor {{\it {Corr}}}
\newcommand \res {{\rm {res}}}
\newcommand \red {{\rm{red}}}
\newcommand \Gal {{\rm {Gal}}}
\newcommand \PGL {{\rm {PGL}}}
\newcommand \Bl {{\rm {Bl}}}
\newcommand \Sing {{\rm {Sing}}}
\newcommand \spn {{\rm {span}}}
\newcommand \Nm {{\rm {Nm}}}
\newcommand \inv {{\rm {inv}}}
\newcommand \codim {{\rm {codim}}}
\newcommand \Div{{\rm{Div}}}
\newcommand \CH{{\rm{CH}}}
\newcommand \sg {{\Sigma }}
\newcommand \DM {{\sf DM}}
\newcommand \Gm {{{\mathbb G}_{\rm m}}}
\newcommand \tame {\rm {tame }}
\newcommand \znak {{\natural }}
\newcommand \lra {\longrightarrow}
\newcommand \hra {\hookrightarrow}
\newcommand \rra {\rightrightarrows}
\newcommand \ord {{\rm {ord}}}
\newcommand \NS {{\rm {NS}}}
\newcommand \Rat {{\mathscr Rat}}
\newcommand \rd {{\rm {red}}}
\newcommand \bSpec {{\bf {Spec}}}
\newcommand \Proj {{\rm {Proj}}}
\newcommand \pdiv {{\rm {div}}}
\newcommand \wt {\widetilde }
\newcommand \ac {\acute }
\newcommand \ch {\check }
\newcommand \ol {\overline }
\newcommand \Th {\Theta}
\newcommand \cAb {{\mathscr A\! }{\it b}}

\newenvironment{pf}{\par\noindent{\em Proof}.}{\hfill\framebox(6,6)
\par\medskip}

\newtheorem{theorem}
[subsection]
{Theorem}
\newtheorem{proposition}[subsection]{Proposition}
\newtheorem{lemma}[subsection]{Lemma}
\newtheorem{corollary}[subsection]{Corollary}
\newtheorem{Remark}[subsection]{Remark}

\theoremstyle{definition}
\newtheorem{definition}[subsection]{Definition}

\title[Selmer group associated to Chow group]{Selmer group associated to the Chow group of certain codimension two cycles }
\author{Kalyan Banerjee, Kalyan Chakraborty}

\address{Kalyan Banerjee @VIT University , Chennai 600127, India.}
\email{kalyan.banerjee@vit.ac.in}
\address{Kalyan Chakraborty @Kerala School of Mathematics, Kozhikode 673571, Kerala, India.}
\email{kalychak@ksom.res.in}


\begin{abstract}
Let $X$ be a surface with geometric genus and irregularity zero which is defined over a number field $K$. Let $\bcX$ denote a smooth spread of $X$ over $\bcO_K[1/f]$ for some element $f\in \bcO_K$ and $A^2$ stands for the group of algebraically trivial cycles on schemes modulo rational equivalence. If $j^*: A^2(\bcX)\to A^2(X)$ be the flat pull-back corresponding to the embedding $j:X\hookrightarrow \bcX$ then we prove that $\im(j^*)(K)/A^2(\bcX)(K)$ is a torsion group. Here $\im(j^*)(K), A^2(\bcX)(K)$ stand for the cycles fixed under the action of the absolute Galois group.
\end{abstract}

\maketitle

\section{Introduction}
Suppose
$X$  be a smooth projective surface defined over a number field $K$ and assume that it can be spread out to a smooth projective scheme $\bcX$ over an affine open subset of the spectrum of the number ring $\bcO_K$. Let $A^2(X)$ denote the group of algebraically trivial cycles of codimension $2$ modulo rational equivalence on $X$.

We recall the definition of algebraic equivalence over $\bcO_K$, which will be used in the sequel.
Let us consider the free abelian group of codimension two cycles on $\bcX$. Two cycles $z_1,z_2$ are said to be algebraically equivalent if there exists a smooth projective curve $C$ defined over $\bcO_K$, two scheme theoretic points $x_0, x_1$ on $C$ and a relative correspondence $\Gamma$ on $C\times_{\bcO_K}\bcX$, such that the intersection
$$
\Gamma.(x_0\times_{\bcO_K}\bcX)-\Gamma.(x_1\times_{\bcO_K}\bcX)=z_1-z_2.
$$
Here $.$ denotes the relative intersection product in the sense of \cite{Fu}.

Let $\bar K$ denote the algebraic closure of $K$ and $G = $ Gal $( {\bar{K}}/K)$ be the absolute Galois group.
Also,  $X_{\bar K}$ denotes the surface
$$
X\times _{K}\bar K
$$
and
$$
\bcX_{\bar K}:=\bcX\times_{\bcO_K}\times \overline{\bcO_K}.
$$
Here $\overline{\bcO_K}$ denotes the integral closure of $\bcO_K$ in $\bar K$.

Let $A^2(\bcX_{\bar K})(K), A^2(X_{\bar K})(K)$ be the $G$-fixed part of the action of $G$ on $A^2(\bcX_{\bar K}),A^2(X_{\bar K})$ respectively,
and $j$ be the embedding of $X_{\bar K}$ into $\bcX_{\bar K}$.
If one considers the flat pullback
$$
j^*: A^2(\bcX_{\bar K})\to A^2(X_{\bar K})
$$
of codimension $2$-cycles, it gives a map
$$
j^*: A^2(\bcX_{\bar K})(K)\to A^2(X_{\bar K})(K).
$$
Then a general question is: what is the cokernel  of this homomorphism?

Mildenhall \cite{Mil} studied this flat pull-back for $X$ to be the self-product of an elliptic curve admitting a complex multiplication. He has shown that the kernel of this flat-pullback over any number field $K$ is finite. We use  Mildenhall's result to derive the information about the quotient $\im(j^*)(K)[n]/A^2(\bcX)(K)[n]$ at the level of $n$-torsions of this homomorphism for the case where $X$ is the self-product of a CM elliptic curve.


Let $E$ be an elliptic curve with complex multiplication by the ring of integers of a number field $K$ and $N$ be  it's discriminant (for this we fix an Weirstrass equation for $E$ once and for all). We denote $E\times E$ by $X$
and suppose that there exists a smooth spread of $X$, say $\bcX$ defined over $\bcO_K[1/6N]$.
Let
$j^*$ be the flat pull-back at the level of $A^2$ induced by the embedding
$$j:X_{\bar K}\to \bcX_{\bar K}.$$
Then the main result is:

{\bf{Theorem \ref{theorem1}}}: The group $S^{6N}(\Sigma_{\bar K})$ is a colimit of the unramified cohomology groups
$$
 H^1(G,\Sigma_{\bar K}[6N])\to H^1(I_v,\Sigma_v[6N])
 $$
 here $G$ is the Galois group $Gal(\bar K/L)$ for a finite extension $L$ of $K$ and $I_v$ is the inertia subgroup of the Galois group $G$ corresponding to a finite place $v$ and hence
$$
\im(j^*[6N])(K)/A^2(\bcX_{\bar K})[6N](K)
$$
is a colimit of unramified cohomology groups. 

\medskip

The main tools used here are the Galois module structure of the Chow group of $X_{\bar K}$ and that of $\bcX_{\bar K}$ and
the Galois cohomology of the groups $A^2(\bcX_{\bar K})$ and that of $A^2(X_{\bar K})$. The proof involves similar techniques as to show that the Selmer group of an abelian variety defined over a number field is finite.

Coombes  \cite{Co}
proved that for a surface $X$ with geometric genus and irregularity zero,  $A^2(X)$
is finite under the assumption that $A^2(X_{\bar K})=0$.

Towards proving  our main result  Theorem \ref{theorem1}, we start with a surface $X$ of geometric genus and irregularity zero defined over $K$  which satisfies the condition that the map:
$$
\Pic(X)\to \Pic(X_{\bar K})\to \NS(X_{\bar K})
$$
is surjective,
$$
H^2(\bcX_{\bcO_K[1/f], \mathfrak{p}}, \bcO_{\bcX_{\bcO_K[1/f]},\mathfrak{p}})=0
$$
for all $\mathfrak p$ in $\Spec(\bcO_K[1/f])$, here $f$ is some element in $\bcO_K$ and
$$X(K)\neq \emptyset\;.$$
Then we prove by using the result of \cite{CTR}[lemma 3.3], that:
\medskip

{\bf{Theorem \ref{theorem2}}}: The group
$$
\im(j^*)(K)/A^2(\bcX_{\overline{\bcO_K[1/f]}})(K)
$$
is a torsion group.

\medskip

Here we do not assume the vanishing of $A^2(X_{\bar K})$ as in \cite{Co}. The above theorem is important to prove the finiteness or triviality of $A^2(X_{\bar K})(K)$. It says that, atleast to prove that $\im(j^*)(K)\otimes_{\ZZ} \QQ=\{0\}$ it is enough to prove that
$$
A^2(\bcX_{\overline{\bcO_K[1/f]}})(K)\otimes_{\ZZ} \QQ=\{0\}\;.
$$
This gives some information on how to prove $A^2(X)$ is finite.

\smallskip

{\small \textbf{Acknowledgements:} The authors thank department of atomic energy (DAE) for funding this project and for the hospitality of Harish-Chandra Research Institute, India, where the work has been done. The first author also thanks VIT University Chennai for hosting this project.}

\section{Proof of the theorems}
Let $E$ be as before
having complex multiplication by $\bcO_K$ and $X = E\times E$. Let us fix a Weirstrass equation for the elliptic curve once and for all. Let us consider a spread $\bcE_{\bcO_K[1/6N]}$ of $E$ over $\bcO_K[1/6N]$ which is smooth and denote
$$\bcE_{\bcO_K[1/6N]}\times \bcE_{\bcO_K[1/6N]}$$
by $\bcX$. Then we consider the restriction homomorphism from $A^2(\bcX)\to A^2(X)$. It is known due to Mildenhall's result that the kernel of this restriction map is finite and we name it $\Sigma_K$. Then we have the exact sequence
$$
0\to \Sigma_K \to A^2(\bcX)\to A^2(X)\;.
$$
Now consider the sequence at the level of $\bar{K}$, namely
$$
0\to \Sigma_{\bar K}\to A^2(\bcX_{\bar K})\to A^2(X_{\bar K})\;.
$$
Note that $G$
acts naturally on each member of the above short exact sequence. Also
if the inclusion of $X_{\bar K}\hookrightarrow \bcX_{\bar K}$ be denoted by $j$, then the pullback map $j^*$ is from $A^2(\bcX_{\bar K})$ to $A^2(X_{\bar K})$.

Therefore one has the natural long exact sequence on the group cohomology level of $G$ for these Galois modules,
$$
0\to \Sigma_{\bar K}^G\to A^2(\bcX_{\bar K})^G\to  \im(j^*)^G\to H^1(G,\Sigma_{\bar K})\to H^1(G,A^2(\bcX_{\bar K}))\to H^1(G,\im(j^*))\;.
$$
Here $M^G$, for a $G$-module $M$, denotes the group of $G$-invariants in $M$, i.e.
$$
\{m\in M|g.m=m, \forall g\in G\}\;.
$$
Let us denote the groups
$$
A^2(X_{\bar K})^G, A^2(\bcX_{\bar K})^G, \im(j^*)^G
$$
as
$$
A^2(X_{\bar K})(K), A^2(\bcX_{\bar K})(K), \im(j^*)(K)\;
$$
respectively.
Also for notational convenience we continue to denote $\Sigma_{\bar K}^G$ as $\Sigma_K$.
Then we have the following exact sequence
$$
0\to \im(j^*)(K)/A^2(\bcX_{\bar K})(K)\to H^1(G,\Sigma_{\bar K})\to H^1(G,A^2(\bcX_{\bar K}))\to H^1(G,\im(j^*))\;.
$$
Let $v$ be a place of $K$ and
$K_v$ be
the completion of $K$ at $v$.
Let $\bar K_v$  be the algebraic closure of $K_v$ and we embed $\bar K$ into $\bar K_v$. This embedding gives an injection of
 $\Gal(\bar K_v/K_v)=G_v$ into $\Gal(\bar K/K)=G$ and consequently
a homomorphism (considering the Galois cohomology)
$$
H^1(G,\Sigma_{\bar K})\to H^1(G_v,A^2(\bcX_{\overline { \bcO_{K_v}[1/6N]}}))\;.
$$
Again for notational convenience we write $A^2(\bcX_{\overline { \bcO_{K_v}[1/6N]}})$ in the above as $A^2(\bcX_{\bar K_v})$.
Then we have the following commutative diagrams:
$$
  \xymatrix{
  \im(j^*)(K)/A^2(\bcX_{\bar K})(K)\ar[r]^-{} \ar[dd]^-{} & H^1(G,\Sigma_{\bar K}) \ar[r]^-{} \ar[dd]_-{}
  &   H^1(G,A^2(\bcX_{\bar K})) \ar[r]^-{} \ar[dd]_-{}
  & H^1(G,\im(j^*))  \ar[dd]_-{}  \
  \\ \\
 \im(j^*)(K_v)/A^2(\bcX_{\bar {K_v}})(K_v)\ar[r]^-{} & H^1(G_v,\Sigma_{\bar K_v}) \ar[r]^-{}
    & H^1(G_v,A^2(\bcX_{\bar K_v})) \ar[r]^-{}
  & H^1(G_v,\im(j_v^*))
  }
$$
Let us now focus on
$$
H^1(G,\Sigma_{\bar K})\to \prod_v H^1(G_v,A^2(\bcX_{\bar K_v}))\;
$$
and consider the sequence of $n$-torsion subgroups of $\Sigma_{\bar K}, A^2(\bcX_{\bar K}), A^2(X_{\bar K})$ given by:
$$
0\to \Sigma_{\bar K}[n]\to A^2(\bcX_{\bar K})[n]\to A^2(X_{\bar K})[n]\;.
$$
Here $A[n]$ for an abelian group $A$, denotes the group of $n$-torsions of $A$.
Then considering the above groups as $G$-modules we have a homomorphism at the level of Galois cohomology  given by:
$$
H^1(G,\Sigma_{\bar K}[n])\to \prod_v H^1(G_v,A^2(\bcX_{\bar K_v}))[n]\;.
$$
\begin{definition}
The kernel of this map is defined to be the $n$-Selmer group associated to the restriction homomorphism $A^2(\bcX_{\bar K})\to A^2(X_{\bar K})$, at the level of $n$-torsions in the  group of algebraically trivial  codimension $2$-cycles and it is denoted by $S^n(\Sigma_{\bar K})$.
\end{definition}
Let $Alb(X_{\bar K})$ be the Albanese variety such that there exists a natural (universal) homomorphism of abelian groups from $A^2(X_{\bar K})$ to $Alb(X_{\bar K})$. Now $X_{\bar K}$ is an abelian variety $Alb(X_{\bar K})\cong X_{\bar K}$. Since the following argument is more general in nature, that is, it works for any smooth projective $X_{\bar K}$ and for its albanese variety $Alb(X_{\bar K})$, provided the kernel of
$$
j^*: A^2(\bcX_{\bar K})[n](K)\to A^2(X_{\bar K})[n](K)
$$
is finite, we do not use the isomorphism $Alb(X_{\bar K})\cong X_{\bar K}$. Specifically this is required to prove the analogous result as stated in Remark \ref{remark1}.

Let's consider the commutative diagram:
$$
  \diagram
  H^1(G,Alb(X_{\bar K})[n])\ar[dd]_-{} \ar[rr]^-{} & & \prod_v H^1(G_v,Alb(X_{\bar K_v}))[n] \ar[dd]^-{} \\ \\
 H^1(G,A^2(X_{\bar K})[n])\ar[rr]^-{} & & \prod_v H^1(G_v,A^2(X_{\bar K_v}))[n]
  \enddiagram
$$
Now by Roitman's theorem \cite{R2}, the groups $Alb(X_{\bar K})[n]$ and $A^2(X_{\bar K})[n]$ are isomorphic as Galois modules and therefore the group cohomologies are isomorphic. Thus the left vertical arrow in the above diagram is an isomorphism. Let

$$S^n(Alb(X_{\bar K})/K):=\ker(H^1(G_K,Alb(X_{\bar K})[n])\to \prod_v H^1(G_{K_v},Alb(X_{\bar K_v})[n]))$$
here $v$ varies over all finite places of $K$. Similarly 
$$S^n(A^2(X_{\bar K})/K):=\ker(H^1(G_K,A^2(X_{\bar K})[n])\to \prod_v H^1(G_{K_v},A^2(X_{\bar K_v})[n]))\;.$$
 If we take an element
 in
$S^n(Alb(X_{\bar K})/K)$, then by the commutativity of the above diagram, the image of the element under the left vertical homomorphism is in $S^n(A^2(X_{\bar K})/K)$. Now we prove our main result which has already been stated in the introduction.

\begin{theorem}
\label{theorem1}
The group $S^{6N}(\Sigma_{\bar K})$ is a colimit of the unramified cohomology groups
$$
 H^1(G,\Sigma_{\bar K}[6N])\to H^1(I_v,\Sigma_v[6N])
 $$
 here $G$ is the Galois group $Gal(\bar K/L)$ for a finite extension $L$ of $K$ and $I_v$ is the inertia subgroup of the Galois group $G$ corresponding to a finite place $v$ and hence
$$
\im(j^*[6N])(K)/A^2(\bcX_{\bar K})[6N](K)
$$
is a colimit of unramified cohomology groups. 
\end{theorem}
\begin{proof}
Let $n$ be a positive integer.  Let us
consider the diagram
$$
  \diagram
  H^1(G,\Sigma_{\bar K}[n])\ar[dd]_-{} \ar[rr]^-{} & & \prod_v H^1(G_v,\Sigma_{\bar K_v})[n] \ar[dd]^-{} \\ \\
 H^1(G,A^2(\bcX_{\bar K})[n])\ar[rr]^-{} & & \prod_v H^1(G_v,A^2(\bcX_{\bar K_v}))[n]
  \enddiagram
 $$
Suppose that some element is there in $S^n(\Sigma_{\bar K}/K)$.

Consider the following commutative squares:
$$
  \diagram
  \Sigma_{K_v}/n\Sigma_{K_v}\ar[dd]_-{} \ar[rr]^-{} & & H^1(G_v,\Sigma_{\bar K_v}[n]) \ar[dd]^-{} \\ \\
  A^2(\bcX_{\bar K_v})(K_v)/nA^2(\bcX_{K_v})(K_v) \ar[rr]^-{} & & H^1(G_v,A^2(\bcX_{\bar K_v})[n])
  \enddiagram
$$
Let $\eta$ be in the kernel of
 $$
 H^1(G,\Sigma_{\bar K}[n])\to H^1(G_v,A^2(\bcX_{\bar K_v}))[n]\;.
 $$
 Then it follows by the exactness of the sequence, induced by the long exact sequence corresponding to the short exact sequence of Galois modules given as : $$0\to \Sigma_{\bar K}[n]\to A^2(\bcX_{\bar K_v})\stackrel{n}{\to} A^2(\bcX_{\bar K_v})\to 0 $$
 $$
 A^2(\bcX_{K_v})(K_v)/nA^2(\bcX_{K_v})(K_v) \to H^1(G_v,A^2(\bcX_{\bar K_v})[n])\to H^1(G_v,A^2(\bcX_{\bar K_v}))[n]
 $$
 that there exists an element $z$ in $A^2(\bcX_{K_v})(K_v)$ such that
 $$
 \phi(\eta)(\sigma)=\sigma.z-z
 $$
 for all $\sigma$ in $G_v$. Here
 $\phi$ is from $Z^1(G,\Sigma_{\bar K}[n])$ to $Z^1(G_v,A^2(\bcX_{\bar K_v})[n])$ and $\eta$ is a co-cycle such that it's cohomology class is in the kernel the homomorphism induced by $\phi$ at the level of cohomology. For simplicity as before we denote $\bcX_{\bar K_v},\Sigma_{\bar K_v}$ by $\bcX_v, \Sigma_v$ respectively. In particular for all $\sigma$ in the inertia group $I_v$, we have
 $$
 \phi(\eta)(\sigma)=\sigma.z-z\;.
 $$
Let $v$ be a finite place such that $v$ does not divide $n$ and $Alb(X_{\bar K}), X_{\bar K}$ have good reduction at $v$. We consider the specialization homomorphism from $A^2(X_v)$ to $A^2(X_v')$, where $X_v'$ is the reduction of $\bcX_v$ at $v$. Then it follows that the image of
 $$
 \sigma.z-z
 $$
in $A^2(X_v)$ goes to zero  under the specialization homomorphism for all $\sigma$ in $I_v$. But on the other hand
$$
\sigma.z-z
$$
is an $n$-torsion for each $\sigma$ in $G_v$ (as $\eta$ is an $n$-torsion), so by Roitman's theorem on {\it torsion}, the image of the element $\sigma.z-z$ in $A^2(X_v)$ corresponds to an $n$-torsion on $Alb(X_v)$. By the previous argument
 this $n$-torsion on $Alb(X_v)$ is mapped to zero under $Alb(X_v)\to Alb(X_v')$. But we know that the $n$-torsions of $Alb(X_v)$ are embedded in $Alb(X_v')$ (for $v$ which does not divide $n$, this follows from the theory of formal groups over $v$-adic numbers). Therefore this $n$-torsion on $A^2(X_v)$ is zero (this is because of the injectivity of the albanese map on $n$-torsions) and consequently
$$
 \sigma.z-z\in \Sigma_{v}[n]
$$
 for all $\sigma\in I_v$ and
 $$\sigma.z-z=0$$
 in $\Sigma_v[n]$, for all $\sigma\in I_v$ (where $v$ does not divide $n$).

This implies  $S^n(\Sigma_{\bar K})$ consists of elements which are unramified for all but finitely many places $v$, i.e., the image of the elements in $S^n(\Sigma_{\bar K})$ under the map
 $$
 H^1(G,\Sigma_{\bar K}[n])\to H^1(I_v,\Sigma_v[n])
 $$
 is zero for all but finitely many places $v$.

 Hence the following variance of lemma \cite{Sil}[lemma 4.3, chapter X] tells us that $S^{6N}(\Sigma_{\bar K})$ is a colimit of this unramified cohomology groups..
\begin{lemma}\cite{Sil}[Lemma 4.3, Chapter 10]
Let $L$ be a finite extension of the number field $K$. Let $M$ be the  finite $G=\Gal(\bar K/L)$ module $\Sigma_L$ and $S$ be a set of finitely many places in $L$. Consider
 $$
 H^1(G,M;S)
 $$
consisting of all elements $\eta$ in $H^1(G,M)$, which are unramified outside $S$, that is in the kernel
$$\ker(H^1(G,M)\to \prod_{v\not\in S} H^1(I_v,M_v))$$
here $M_v=\Sigma_v$, for a place in $L$. Then $H^1(G,M;S)$ is finite.
\end{lemma}
Thus this result, applied to $M=\Sigma_{L}[6N]$ for any finite extension $L$ over $K$ (is finite by  Theorem 1.1 in \cite{CTR}), we have that $S^{6N}(\Sigma_L)$ is finite. Also observe that
$$
\Sigma_{\bar K}=\cup_{K\subset L} \Sigma_L\;,
$$
$$
G=\cup_{K\subset L} \Gal(\bar K/L)\;.
$$
Therefore $H^1(G,\Sigma_{\bar K}[6N])$ is isomorphic to the colimit of
$$
H^1(\Gal(\bar K/L),\Sigma_L[6N]).
$$
Hence $S^{6N}(\Sigma_{\bar K})$ is isomorphic to the colimit of $S^{6N}(\Sigma_L)$, all of which are finite unramified cohomology groups. Consequently $$\im(j^*[6N])(K)/A^2(\bcX)[6N](K)$$ is a colimit of unramified cohomology groups. The proof actually follows from the finiteness of $S^{6N}(\Sigma_L)$.
\end{proof}
\begin{Remark}
In the previous theorem \ref{theorem1}, it is interesting to see whether the groups $S^{6N}(\Sigma_L)$ are subgroups of the group $S^{6N}(\Sigma_{\bar K})$.
\end{Remark}
\begin{Remark}\label{remark1}
The analogue of
Mildenhall's result was proved for a Fermat quartic surface in \cite{Ot}. Thus for the Fermat quartic surface too,  Theorem \ref{theorem1} is true.
\end{Remark}
Now by lemma 3.3 in \cite{CTR}, if
$\bcO_L[1/f]$ is such that $\bcX_L:=\bcX_{\bcO_L[1/f]}$ is smooth and the following conditions are true:
\begin{eqnarray*}
H^2(\bcX_{L,\mathfrak{p}},\bcO_{\bcX_{L,\mathfrak{p}}})&=&0 ~~\mbox{for all} ~\mathfrak{p}\in \Spec(\bcO_L[1/f]),\\
X(K)&\neq& \emptyset,\\
\Pic(X_L)\to \Pic( X_{\bar K})&\to& \NS( X_{\bar K}) ~~\mbox{is surjective},
\end{eqnarray*}
 then for a finite extension $L$ of $K$, $A^2(\bcX_L)\to A^2(X_L)$ has finite kernel. This leads us to prove:
\begin{theorem}
\label{theorem2}
Under the above conditions
$$
\im(j^*)(K)/A^2(\bcX_{\bar K})(K)
$$
is torsion and is described by a colimit of unramified cohomology groups

$$
 H^1(G,\Sigma_{\bar K}[n])\to H^1(I_v,\Sigma_v[n])\;.
 $$
 
\end{theorem}
\begin{proof}
The proof goes verbatim
as  Theorem \ref{theorem1}.
\end{proof}
\begin{Remark}
For a surface with geometric genus and irregularity zero, such that the above conditions as in Theorem \ref{theorem1} are satisfied, one has that
$$
\im(j^*)(K)/A^2(\bcX_{\bar K})(K)
$$
is torsion. Therefore tensoring with $\QQ$ gives
$$
A^2(\bcX_{\bar K})(K)\otimes \QQ\to \im(j^*)(K)\otimes \QQ
$$
is surjective. Therefore to prove the triviality  for $\im(j^*)(K)\otimes \QQ$, it is enough to prove that
$$
A^2(\bcX_{\bar K})(K)\otimes \QQ
$$
is trivial. This gives some information about the structure of of $A^2(X_{\bar K})$.
\end{Remark}

\end{document}